\documentclass[psamsfonts]{amsart}
 
\font\tenbb=msbm10 scaled 1000 
 
\font\sevenbb=msbm7 scaled 1000 
 
\font\fivebb=msbm5 scaled 1000 
 
\newfam\bbfam 
 
\textfont\bbfam=\tenbb 
 
\scriptfont\bbfam=\sevenbb 
 
\scriptscriptfont\bbfam=\fivebb 
 
\def\bb{\fam\bbfam\tenbb} 
  
\input{psfig.tex} 
\input{diagrams.sty} 
\input{psfig.sty} 
 
\usepackage{latexsym,amsfonts,amssymb}
\usepackage{amsbsy} 
 \usepackage{amsfonts} 
 \usepackage{amsmath} 
 \newtheorem{thm}{Theorem}[section] 
 \newtheorem{prop} {Proposition}[section] 
  
 \newtheorem{lem}{Lemma}[section] 
 \newtheorem{df} {Definition}[section] 
 \newtheorem{rem}{Remark}[section]  
   
\def\z{{\bb Z}} 
  
 \newcommand{\edim}{\hspace*{\fill}$\square$ \vspace{10pt}} 
  
\def\rea{{\bb R}}

 \newcommand{\pre}{\textit{Proof: }}  
  
 \newcommand{\qw}{\hspace{-10pt}}  
 \newcommand{\vsp}{\vspace{10pt}} 
\newcommand{\sg}{\sigma}

\newcommand{\sjp}{\sigma_{j}}

\newcommand{\sip}{\sigma_{i}}
\newcommand{\siip}{\sigma_{i+1}}

\newcommand{\tre}{\sigma_1 \sg_2 \cdots \sigma_{n-1}^2 \cdots  \sg_2
  \sigma_1}

\title{On presentations of  surface braid groups}

\author{Paolo Bellingeri}

\address{Institut Fourier, BP 74,  
Univ.Grenoble I, Math\'ematiques, 38402 Saint-Martin-d'H\`eres cedex, France}  

\email{paolo.bellingeri@ujf-grenoble.fr}
\thanks{Partially supported by the MENRT grant.\\
\mbox{\hspace{11pt}}{\it Mathematics Subject Classification:}
Primary: 20F36. Secondary: 57N05.\\
\mbox{\hspace{11pt}}{\it Key words}: Braid, Surface, Presentation}

\begin{document}










\maketitle

\begin{abstract} 
We give presentations of  braid
groups and pure braid groups on  surfaces
and we show some properties of surface pure braid groups.
\end{abstract}



\section{Presentations for  surface braids}

Let $F$ be an orientable surface and let
$\mathcal{P}=\{P_1, \dots, P_n\}$ be  
a set of $n$ distinct points of $F$. 
A \emph{geometric braid} on $F$ based at 
$\mathcal{P}$ is an $n$-tuple 
$\Psi=(\psi_1, \dots, \psi_n)$ of paths
$\psi_i: [0, 1] \to F$ such that
\begin{itemize}
\item $\psi_i(0)= P_i, \;  i=1 \dots, n;$
\item $\psi_i(1) \in \mathcal{P},  \;  i=1 \dots, n;$
\item $\psi_1(t), \dots, \psi_n(t)$ are  distinct points of
  $F$ for all $t \in [0, 1]$.
\end{itemize}
The usual  product of paths defines a group structure
on the set of  braids up  to homotopies among braids.
This group, denoted  $B(n,F)$, does not depend on the choice of  $\mathcal{P}$
and it is called  the braid group on $n$ strings on $F$.
On the other hand, let be $F_nF=F^n \setminus \Delta$, where $\Delta$ is the big diagonal, i.e. 
the $n$-tuples $x=(x_1, \dots x_n)$ for which $x_i=x_j$ for some $i \not= j$.
There is a natural  action of $\Sigma_n$ on $F_nF$ by permuting coordinates. We call the orbit space
  $\hat{F}_nF= F_nF/\Sigma_n$
{\em configuration space}. Then the braid group
$B(n,F)$ is isomorphic to $\pi_1(\hat{F}_nF)$.
We recall that the pure braid group $P(n,F)$ 
on $n$ strings on $F$ is the kernel of the natural projection 
of $B(n, F)$ in the permutation group $\Sigma_n$. This group is isomorphic to 
$\pi_1(F_nF)$. 
The first aim of this article is to give (new) presentations for braid 
groups on orientable surfaces. 

A $p$-punctured surface of genus $g\ge 1$ is the surface obtained
by deleting $p$ points on a closed  surface of genus $g\ge 1$.

\begin{thm} \label{my:th}
Let   $F$  be an 
orientable $p$-punctured surface of genus $g\ge 1$. The group  $B(n, F)$ admits
the following presentation (see also section \ref{gen:rel}):  

\begin{itemize}
\item Generators:
$\sg_1,\dots,\sg_{n-1}, a_1, \dots,a_{g}, b_1, \dots,b_{g}, z_1,\dots, z_{p-1} \,.$
\item Relations:
\begin{itemize}
\item  Braid relations, i.e.
\begin{eqnarray*}
 \sip \siip \sip &=&  \siip \sip \siip \, ;\\
\sip \sjp &=& \sjp \sip \quad   \mbox{for} \;| i-j | \ge 2  \, .
\end{eqnarray*}
\item Mixed relations:
\end{itemize}
\end{itemize}
\begin{eqnarray*}
&(R1) &\;a_r \sg_i=\sg_i  a_r      \quad  (1 \le r \le g;\; i\not= 1 )\, ; \\
&     &\;b_r \sg_i=\sg_i  b_r      \quad  (1 \le r \le g;\; i\not= 1 )\, ; \\
&(R2) &\;\sigma_1^{-1} a_r \sigma_1^{-1} a_{r}= a_r \sigma_1^{-1} a_{r}\sigma_1^{-1} \quad  (1 \le r \le g) \, ; \\
&     &\;\sigma_1^{-1} b_r \sigma_1^{-1} b_{r}= b_r \sigma_1^{-1} b_{r}\sigma_1^{-1} \quad  (1 \le r \le g) \, ; \\
&(R3)&\;  \sigma_1^{-1} a_{s} \sigma_1 a_r = a_r \sigma_1^{-1} a_{s} \sigma_1\quad (  s < r) \, ; \\
&    &\;  \sigma_1^{-1} b_{s} \sigma_1 b_r = b_r \sigma_1^{-1} b_{s} \sigma_1\quad (  s < r) \, ; \\
&    &\;  \sigma_1^{-1} a_{s} \sigma_1 b_r = b_r \sigma_1^{-1} a_{s} \sigma_1\quad (  s < r) \, ; \\
&    &\;  \sigma_1^{-1} b_{s} \sigma_1 a_r = a_r \sigma_1^{-1} b_{s} \sigma_1\quad (  s < r) \, ; \\
&(R4) &\; \sigma_1^{-1} a_r \sigma_1^{-1} b_{r} = b_{r} \sigma_1^{-1} a_r \sigma_1 \quad(1 \le r \le g)\, ;\\ 
&(R5) &\; z_j\sg_i= \sg_i z_j        \quad  (i \not=n-1, j=1, \dots, p-1 )\, ;\\ 
&(R6)&\;  \sigma_1^{-1} z_i \sigma_1 a_r =a_r\sigma_1^{-1} z_i\sigma_1    \quad  (1 \le r \le g;\; i=1,
\dots, p-1; \; n>1 )\,; \\
&    &\;  \sigma_1^{-1} z_i \sigma_1 b_r =b_r\sigma_1^{-1} z_i\sigma_1    \quad  (1 \le r \le g;\; i=1,
\dots, p-1; \; n>1 )\,; \\
&(R7)&\;\sg_{1}^{-1} z_j \sg_{1} z_l=z_l\sg_{1}^{-1} z_j \sg_{1} \quad (j=1, \dots, p-1 , \; j<l )\, ;\\
&(R8)&\;\sg_{1}^{-1} z_j \sg_{1}^{-1} z_j=z_j\sg_{1}^{-1}  z_j \sg_{1}^{-1} \quad (j=1, \dots, p-1)  \,    .
\end{eqnarray*}
\end{thm}

\begin{thm}\label{men}
Let  $F$ be  a closed orientable surface of genus $g\ge1$. The group $B(n,F)$ admits the following presentation: 
\begin{itemize}
\item Generators:
$\sg_1,\dots,\sg_{n-1}, a_1, \dots,a_{g} \,  b_1, \dots, b_{g}.$
\item Relations:
\begin{itemize}
\item Braid relations as in Theorem \ref{my:th}.

\item Mixed relations:
\end{itemize}
\end{itemize}
\begin{eqnarray*}
&(R1) &\;a_r \sg_i=\sg_i  a_r      \quad  (1 \le r \le g;\; i\not= 1 )\, ; \\
&     &\;b_r \sg_i=\sg_i  b_r      \quad  (1 \le r \le g;\; i\not= 1 )\, ; \\
&(R2) &\;\sigma_1^{-1} a_r \sigma_1^{-1} a_{r}= a_r \sigma_1^{-1} a_{r}\sigma_1^{-1} \quad  (1 \le r \le g) \, ; \\
&     &\;\sigma_1^{-1} b_r \sigma_1^{-1} b_{r}= b_r \sigma_1^{-1} b_{r}\sigma_1^{-1} \quad  (1 \le r \le g) \, ; \\
&(R3) &\;  \sigma_1^{-1} a_{s} \sigma_1 a_r = a_r \sigma_1^{-1} a_{s} \sigma_1 \quad (  s < r) \, ; \\
&     &\;  \sigma_1^{-1} b_{s} \sigma_1 b_r = b_r \sigma_1^{-1} b_{s} \sigma_1 \quad (  s < r) \, ; \\
&     &\;  \sigma_1^{-1} a_{s} \sigma_1 b_r = b_r \sigma_1^{-1} a_{s} \sigma_1 \quad (  s < r) \, ; \\
&     &\;  \sigma_1^{-1} b_{s} \sigma_1 a_r = a_r \sigma_1^{-1} b_{s} \sigma_1 \quad (  s < r) \, ; \\
&(R4) &\;  \sigma_1^{-1} a_r \sigma_1^{-1} b_{r} = b_{r} \sigma_1^{-1} a_r \sigma_1 \quad(1 \le r \le g) \,     ; \\
&(TR) &\; [a_1,b_1^{-1}] \cdots [a_{g},b_{g}^{-1}]=\tre \, ,
\end{eqnarray*}
where $[a, b]:=a b a^{-1} b^{-1}$.
\end{thm}

\vsp

We may assume that Theorem \ref{my:th} provides also a presentation for  $B(n,F)$, when $F$ an orientable surface with $p$ 
boundary components. 
When $F$ is a closed orientable surface, our presentations
are similar to Gonz\'alez-Meneses' presentations, but the number of  relations is smaller.
We recall also that  the first 
presentations of braid groups on closed surfaces
were found  by Scott (\cite{sco}), afterwards revised by Kulikov and Shimada
(\cite{kul}). At our knowledge, the case of punctured surfaces is new in the literature. 
Our proof is  inspired by
Morita's combinatorial proof  for the classical
presentation of Artin's braid group (\cite{mor}). 
We will explain this approach while proving Theorem \ref{my:th}. 
After that we will show how to make  this technique  fit for obtaining
Theorem \ref{men}.
The last part of the article concerns the study of 
surface pure braids groups, for  $F$  an orientable surface.    
We provide in Theorem \ref{pure} a  homogeneous presentation for $P(n,F)$, very close to the 
standard presentation of the pure braid group $P_n$ on the disk. Several results on surface pure braid groups are deduced.

\qw {\bf Acknowledgments.} The author is grateful to John Guaschi
for his preprint on surface pure braids and to   Barbu Berceanu, Louis Funar, Juan
Gonz\'alez-Meneses,  Stefan Papadima and Vlad Sergiescu for useful discussions and suggestions. 
Part of this work was done 
during the author's visit to the I.M.A.R. of Bucharest,
whose support and hospitality are gratefully acknowledged.

\section{Preliminaries}

\subsection{Fadell-Neuwirth fibrations} \label{section}

The main tool one uses  is the Fadell-Neuwirth fibration, with 
its generalisation and the corresponding exact sequences. 
As observed in \cite{fan},
if $F$ is a  surface (closed or punctured, orientable or not), the map  $\theta: F_{n}F \to F_{n-1}F$ defined by
$$
\theta(x_1, \dots, x_n)=(x_1, \dots, x_{n-1})
$$
is a fibration with fiber $F \setminus \{x_1, \dots, x_{n-1}\}$.
The exact  homotopy  sequence of the fibration gives us the exact sequence
$$
\cdots \pi_2(F_nF) \to \pi_2(F_{n-1}F) \to \pi_1(F \setminus \{x_1, \dots, x_{n-1} \})
$$
$$
\to
P(n,F) \to P(n-1,F) \to 1.
$$ 
Since a punctured surface (with at least one puncture) has the homotopy type of a one dimensional complex, we deduce
$$
\pi_k(F_nF) \cong \pi_k(F_{n-1}F) \cong  \cdots \cong  \pi_k(F)  , \quad k\ge 3
$$
and
$$
\pi_2(F_nF)  \subseteq \pi_2(F_{n-1}F)  \subseteq  \cdots  \subseteq \pi_2(F) \, .
$$
If  $F$ is an orientable surface and $F \not= S^2$, all  higher homotopy groups are trivial. Thus, 
if $F$ is an orientable surface different from the sphere we can conclude that 
there is an exact sequence
$$
(PBS) \quad 1 \rTo \pi_1(F \setminus \{x_1, \dots, x_{n-1}\}) \rTo P(n,F) \rTo^{\theta} P(n-1,F) \to 1,
$$
where $\theta$ is the map that ``forgets'' the last path pointed at $x_n$.

The problem of the existence of a section for $(PBS)$ has been completely solved in \cite{gua}.
It is possible to show that $\theta$ admits a section, when $F$ has punctures.
On the other hand,
when $F$ is a closed orientable surface of genus  $g \ge 2$, (PBS) splits if and only if  $n=2$. 
An explicit section is shown  in \cite{bir2} in the case of the torus.


\subsection{Geometric interpretations of generators and relations} \label{gen:rel}

Let $F$ be an orientable surface. Let 
 $\widetilde{B}(n, F)$  be the group with the presentation given in Theorem
\ref{my:th}
or  Theorem  \ref{men} respectively.
The geometric interpretation for generators of
$\widetilde{B}(n, F)$,  when  $F$ is a
closed surface of genus $g\ge1$ is the same as  in \cite{gm}, except that we represent $F$ as
a polygon $L$ of $4g$ sides with the standard identification of edges (see also section \ref{sec}).
We can consider  braids  as paths on $L$,
which we draw with the usual ``over and under'' information at  the crossing points.
Figure 1 presents
the generators of $\widetilde{B}(n,F)$ realized as braids on $L$.

\begin{figure}[ht]
\hspace{20pt}\psfig{figure=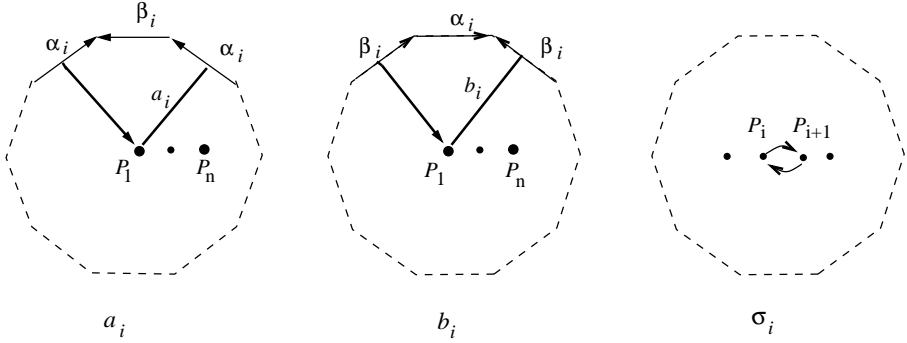,width=12cm} 
 \caption{Generators  as braids (for $F$ an orientable closed surface).}
 \end{figure}

Note  that in the braid  $a_i$ (respectively $b_i$) the only non trivial string is the first one,
which goes  through the the wall $\alpha_i$ (the wall $\beta_i$).
Remark also that $\sg_1 \dots, \sg_{n-1}$ are the 
classical braid generators on the disk.

\begin{figure}[ht] 
 \hspace{20pt}\psfig{figure=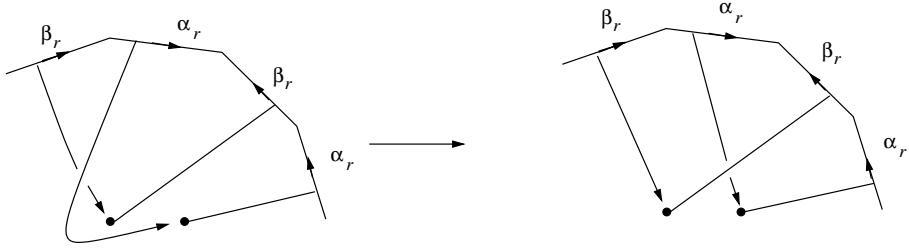,width=12cm} 
 \caption{Geometric interpretation for relation (R4) in Theorem \ref{my:th};
   homotopy between $\sg_1^{-1} a_r \sg_1^{-1} b_{r} $ (on the left) and
 $ b_{r} \sg_1^{-1} a_r \sg_1$ (on the right).}
 \end{figure}

It is easy to check that the relations above hold in $B(n,F)$. 
The non trivial strings of $a_r$ ($b_r$) and $\sg_i$ when $i\not= 1$,
may be considered
to be disjoint
and then  $(R1)$ holds in $B(n,F)$.
On the other hand, $\sg_1^{-1} a_{r} \sg_1^{-1}$ is the braid whose 
 the only non trivial string is the second one,
which goes  through the the wall $\alpha_r$ and   
disjoint from
the corresponding non trivial string of $a_{r}$. Then
$\sg_1^{-1} a_{r} \sg_1^{-1}$ and $a_{r}$ commute. Similarly we have that
$\sg_1^{-1} b_{r} \sg_1^{-1}$ and $b_{r}$ commute and $(R2)$ is verified.
The case of $(R3)$ is similar.
 Figure 2 presents a sketch of a homotopy between
 with $\sg_1^{-1}a_r \sg_1^{-1} b_{r}$ and  $ b_{r} \sg_1^{-1} a_r \sg_1$. Thus, $(R4)$ holds in $B(n,F)$.

Let  $s_r$ (respectively $t_r$) be the first string of $a_r$ (respectively $b_r$), for $r=1, \dots, 2g$, and
consider
all the paths $s_1, t_1, \dots, s_g, t_{g}$.
We cut  $L$ along them and we glue the pieces along the edges of
$L$.
We obtain a new fundamental domain (see Figure 3, for the case of a
surface of genus $2$),
called $L_1$, with vertex $P_1$. On $L_1$ it is clear 
that  $[a_1,b_1^{-1}] \cdots [a_{g},b_{g}^{-1}]$
is equivalent to the braid of Figure 4,  equivalent to the braid
$\sigma_1 \sg_2 \dots \sigma_{n-1}^{2} \dots  \sg_2
\sigma_1$ and then (TR) is verified in $B(n,F)$.

\begin{figure}[ht] 
\hspace{40pt}\psfig{figure=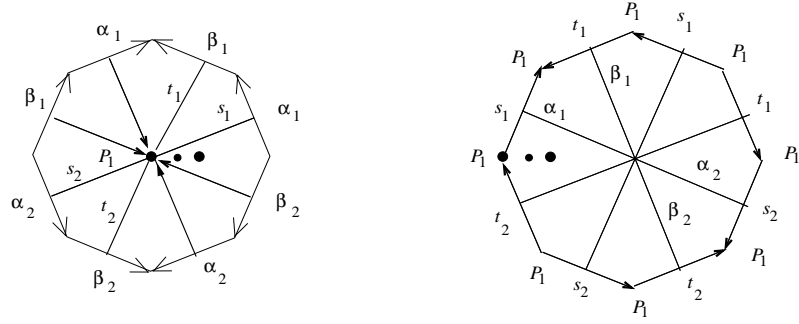,width=10.5cm} 
\caption{The  fundamental domain $L_1$.}
\end{figure}

\begin{figure}[ht] 
\hspace{60pt}\psfig{figure=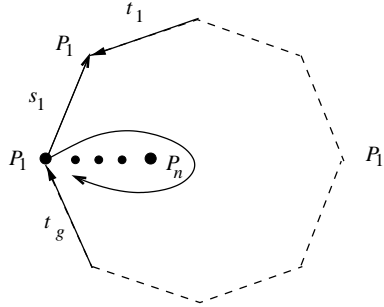,width=5cm} 
\caption{Braid $[a_1,b_1^{-1}] \cdots [a_{g},b_{g}^{-1}]$.}
\end{figure}

There is an analogous  geometric interpretation of generators of
 $\widetilde{B}(n, F)$, for  $F$  an orientable  $p$-punctured surface. 
 The definition of generators $\sg_i, a_j, b_j$ is the same as above.
We only have to add 
 generators $z_i$, where the only non trivial string is the first one,
which is a loop around the
 $i$-th boundary component
(Figure 5), except the $p$-th one.

\begin{figure}[ht] 
 \hspace{20pt}\psfig{figure=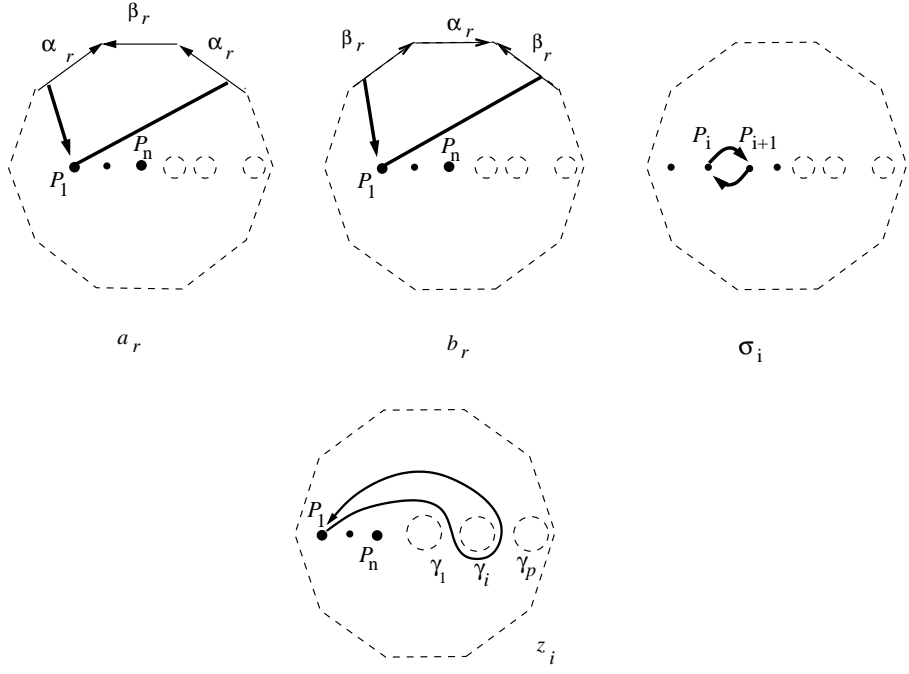,width=12cm} 
 \caption{Generators   as braids (for $F$ an orientable surface
with $p$ boundary components).}
 \end{figure}

As above, relations can be easily checked on corresponding paths (Figure 6).

\begin{figure}[ht] 
 \hspace{50pt}\psfig{figure=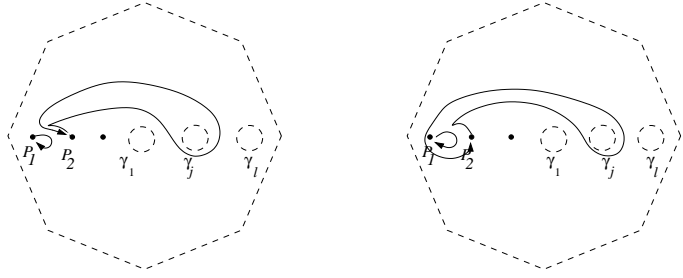,width=9cm} 
 \caption{The braids $\sg_{1}^{-1} z_j \sg_{1}$ and $\sg_{1}^{-1} z_j \sg_{1}^{-1}$.
The non trivial string of  $\sg_{1}^{-1} z_j \sg_{1}$
can be considered  disjoint  from the non trivial string of $z_l$, for $j<l$.
Similarly  the braid  $\sg_{1}^{-1} z_j \sg_{1}^{-1}$ commutes 
with $ z_j$.}
 \end{figure}

Remark that a loop of the first string around the
 $p$-th boundary component can be represented by the geometric braid corresponding to the element
$$
 [a_1,b_1^{-1}] \cdots [a_{g},b_{g}^{-1}]
\sg_1^{-1} \cdots \sg_{n-1}^{-1}\cdots\sg_1^{-1}  z_1  \cdots z_{p-1} \, .
$$

Therefore, one has  natural morphisms $\phi_n: \widetilde{B}(n, F) \to B(n, F)$. 
One further shows that $\phi_n$ are
actually isomorphisms.


\section{Outline of the proof of Theorem \ref{my:th}}


\subsection{The inductive assertion}
We outline the ideas of the proof for $F$ a surface of genus $g$ with one puncture.
One applies an induction on the number $n$ of strands.
For $n=1$, $\widetilde{B}(1, F)=\pi_1(F)=  B(1, F)$,
then $\phi_1$ is an isomorphism.

Consider  the subgroup $B^0(n, F)=\pi^{-1}(\Sigma_{n-1})$ and the  map
$$
\theta: B^0(n, F) \to B(n-1, F)
$$
which ``forgets'' the last string.
Now, let $\widetilde{B}^0(n, F)$  be the subgroup of $\widetilde{B}(n, F)$
generated  
by $a_1, \dots, a_{g}, b_1, \dots, b_{g}, \sg_1, \dots, \sg_{n-2},$ $ \tau_1, \dots , \tau_{n-1},
\omega_1, \dots , \omega_{2g},$
where
\begin{eqnarray*}
\tau_j&=&\sg_{n-1} \cdots \sg_{j+1} \sg_j^2  \sg_{j+1}^{-1} \cdots \sg_{n-1}^{-1} \qquad
(\tau_{n-1}=\sg_{n-1}^2) \, ;\\
\omega_{2r-1}&=& \sg_{n-1}^{-1} \cdots \sg_1^{-1} a_{r} \sg_1
\cdots \sg_{n-1}  \quad r=1, \dots, g \, ;\\
\omega_{2r}&=& \sg_{n-1}^{-1} \cdots \sg_1^{-1} b_{r} \sg_1 
\cdots \sg_{n-1}  \quad  r=1, \dots, g  \, .
\end{eqnarray*}

We construct the following  diagram:
\begin{diagram}
   \widetilde{B}^0(n, F) & \rTo^{\tilde{\theta}} & \widetilde{B}(n-1, F)\\
   \dTo^{\phi_{n|\widetilde{B}^0(n, F)}} &            & \dTo_{\phi_{n-1}} \\
   B^0(n, F) &       \rTo^{\theta} &  B(n-1, F) 
   \end{diagram} 
The map $\tilde{\theta}$ is defined as
$\phi_{n-1}^{-1} \theta
\phi_{n|B^0(n, F)}$. It is well defined, since $\phi_{n-1}$ is an
isomorphism by the inductive assumption, and it is onto.  In fact, 
$\tilde{\theta}(a_i)=a_i, \tilde{\theta}(b_i)=b_i$ for $i=1, \dots, g$
and $\tilde{\theta} (\sg_j)=\sg_j$ for 
$j=1, \dots, n-2$.


\subsection{The existence of a section}

 The morphism $\tilde{\theta}$ has got a natural section
$s: \widetilde{B}(n-1, F) \to \widetilde{B}^0(n, F)$ defined as:
$s(\sg_j)=\sg_j, s(a_i)=a_i, s(b_i)=b_i$
for $j=1, \dots, n-2$ and $i=1, \dots, 2g$.

\begin{rem}
Geometrically 
this section consists of adding a straight strand just to the left  of the puncture.
Generators are sent in corresponding generators.
\end{rem}

Given a group $G$ and a subset $\mathcal{G}$ of elements of 
$G$ we set $\langle \mathcal{G} \rangle$ for the subgroup of $G$ generated  
by $\mathcal{G}$ and $\langle \langle \mathcal{G} \rangle \rangle$
for the subgroup of $G$ normally generated by $\mathcal{G}$.
From now on, given $a, b$ two elements of a group $G$, we set
$a^b= b^{-1} a b$ and ${}^b a =b a  b^{-1}$.

\begin{lem} \label{sect1:bord}
Let  $\mathcal{G}=\{ \tau_1, \dots , \tau_{n-1},
\omega_1, \dots , \omega_{2g} \}$.   
Then $ Ker(\tilde{\theta})=\langle \mathcal{G} \rangle$. 
\end{lem}
\pre We set 
$\beta= \tau_1 \cdots \tau_{n-1}$ $=\sg_{n-1}\cdots  \sg_{2} \sg_{1}^2 \sg_{2} \cdots \sg_{n-1}$
and $\gamma= \beta^{-1}  \tau_1 \beta=
$ $=\sg_{n-1}^{-1} \cdots  \sg_{2}^{-1} \sg_{1}^2 \sg_{2} \cdots \sg_{n-1}$.
By construction we have $\langle \mathcal{G} \rangle \subset
Ker(\tilde{\theta})$.

The existence of a section $s$ implies that
$Ker(\tilde{\theta})=\langle \langle  \mathcal{G} \rangle \rangle$.
In fact, suppose that there is such $x \in Ker(\tilde{\theta})$ 
that $x \notin \langle \langle  \mathcal{G} \rangle \rangle$. Thus,
there is  a word
$x' \not= 1$ on generators $a_1, \dots, a_{g}, b_1,  \dots, b_{g}, \sg_1, \dots, \sg_{n-2},$ of $\widetilde{B}^0(n, F)$
such that $\tilde{\theta}(x')= 1$, because all  other generators 
of $\widetilde{B}^0(n, F)$ are in $\langle \mathcal{G} \rangle$. This is false, since
$x'=s (\tilde{\theta}(x'))$.
To prove that $\langle \mathcal{G} \rangle$ is normal, we need to show that $h^g, {}^g h \in \langle \mathcal{G} \rangle$
for all generators $g$  of $\widetilde{B}^0(n, F)$
and for all $h \in \mathcal{G}$. 

\vspace{5pt}

i) Let $g$ be  one of the  classical  braid generators $\sg_j$, $j=1, \dots, n-2$.
It is clear that   $\tau_i^{\sg_j}$ and $
 {}^{\sg_j}\tau_i$ ($i=1, \dots,n-1)$ belong to $ \langle  \tau_1, \dots , \tau_{n-1}   \rangle$, 
since  it is already true in classical braid groups (\cite{mor}, \cite{vlo}). 
On the other hand, 
$\omega_i^{\sg_j}=
 {}^{\sg_j}\omega_i=\omega_i$ ($i=1, \dots,2g)$.

\vspace{5pt}

ii) Let  $g=a_r$ or $g=b_r$ ($r=1, \dots,g$). Commutativity relations imply
$\tau_j^{g}= {}^{g} \tau_j = \tau_j$ ($ j=2,
\dots n-1$).
Note that
\begin{eqnarray*} 
{}^{a_r}\tau_1={}^{\beta \omega_{2r-1}^{-1}}\gamma \quad \mbox{and} \quad
   \tau_1^{a_r}={}^{\tau_1^{-1} \beta \omega_{2r-1}}\gamma \quad \mbox{for}
   \quad r=1, \dots, g \,  ;\\
{}^{b_r}\tau_1={}^{\beta \omega_{2r}^{-1}}\gamma \quad \mbox{and} \quad
   \tau_1^{b_r}={}^{\tau_1^{-1} \beta \omega_{2r}}\gamma \quad \mbox{for}
   \quad r=1, \dots, g \, .
\end{eqnarray*} 
We show only the first equation (the other is similar).
By iterated application of $[a_r, \sg_1 a_r^{-1}\sg_1]=1$ we obtain:
\begin{eqnarray*} 
{}^{a_r}\tau_1=\sg_{n-1} \cdots \sg_2 a_r \sg_1 a_r^{-1} a_r \sg_1
a_r^{-1}\sg_1\sg_1^{-1} \sg_2^{-1}\cdots \sg_{n-1}^{-1}=\\
=\sg_ {n-1}\cdots \sg_2 a_r \sg_1 
 a_r^{-1} \sg_1 a_r^{-1}\sg_1 
 a_r\sg_1^{-1}\sg_2^{-1}\cdots \sg_{n-1}^{-1}=\\
=\sg_{n-1} \cdots \sg_2\sg_1 a_r^{-1}\sg_1\sg_1 a_r\sg_1^{-1}\sg_2^{-1}\cdots \sg_{n-1}^{-1}=
{}^{\beta \omega_{2r-1}^{-1}}\gamma \, .
\end{eqnarray*} 

Set $a_{2, s}=\sg_1^{-1}a_s\sg_1$ for  $s=1, \dots, g$
and respectively  $b_{2, s}=\sg_1^{-1}b_s\sg_1$ for  $s=1, \dots, g$.
In the same way as above we find  that:
\begin{eqnarray*}
&(RC1) & \;(\sg_1^2)^{a_{r}} = {}^{a_{2,r}} (\sg_1^2)\quad  ( r=1, \dots, g) \, ;
\\
&      & \;(\sg_1^2)^{b_{r}} = {}^{b_{2,r}} (\sg_1^2)\quad  ( r=1, \dots, g) \, ;
\\
&(RC2) & \; {}^{a_{r}}(\sg_1^2) =  (\sg_1^2)^{a_{2,r}\sigma_1^{-2}}\quad  (r=1, \dots,g ) \, ;
\\
&      & \; {}^{b_{r}}(\sg_1^2) =  (\sg_1^2)^{b_{2,r}\sigma_1^{-2}}\quad  (r=1, \dots,g ) \, .
\end{eqnarray*}
Now, remark that relations $(R3)$ and $(R4)$ imply the following relations:
\begin{eqnarray*}
&(R3') &\;a_r \sigma_1 a_{s} \sigma_1^{-1} =\sigma_1 a_{s} \sigma_1^{-1} a_r \quad (r < s) \, ; \\
&      &\;b_r \sigma_1 a_{s} \sigma_1^{-1} =\sigma_1 a_{s} \sigma_1^{-1} b_r \quad (r < s) \, ; \\
&      &\;a_r \sigma_1 b_{s} \sigma_1^{-1} =\sigma_1 b_{s} \sigma_1^{-1} a_r \quad (r < s) \, ; \\
&      &\;b_r \sigma_1 b_{s} \sigma_1^{-1} =\sigma_1 b_{s} \sigma_1^{-1} b_r \quad (r < s) \, ; \\
&(R4') &\;b_r \sigma_1^{-1} a_{r-1} \sigma_1=\sigma_1^{-1} a_{r-1} \sigma_1^{-1} a_r \quad(1 \le r \le  g) \,  ;\\
\end{eqnarray*}
Relations $(RC1),  (RC2) , (R3'), (R4')$  combined with relations $(R2),  (R3), (R4)$ give:
\begin{eqnarray*}
& &\;   {}^{a_r}a_{2, s}=a_{2, s}  \quad (s<r) \, ;\\
& &\;   {}^{b_r}a_{2, s}=a_{2, s}  \quad (s<r) \, ;\\
& &\;   {}^{a_r}b_{2, s}=b_{2, s}  \quad (s<r) \, ;\\
& &\;   {}^{b_r}b_{2, s}=b_{2, s}  \quad (s<r) \, ;\\
& &  \; a_{2, r}^{a_r}= \quad {}^{a_{2,r} \sigma_1^{-2}}a_{2,r} \quad (1 \le r \le  g )\; ; \\
& &  \; b_{2, r}^{b_r}= \quad {}^{b_{2,r} \sigma_1^{-2}}b_{2,r} \quad (1 \le r \le  g )\; ; \\
& & \; {}^{a_r}a_{2, r}={}^{\sigma_1^{2}} a_{2, r} \quad  (1 \le r \le g)\; ;\\
& & \; {}^{b_r}b_{2, r}={}^{\sigma_1^{2}} b_{2, r} \quad  (1 \le r \le g)\; ;\\
& & \; a_{2, s}^{a_{r}}={}^{[a_{2,r},
 \sigma_1^{-2}]}(a_{2, s})
\quad(r<s ) \;; \\
& & \; b_{2, s}^{a_{r}}={}^{[a_{2,r},
 \sigma_1^{-2}]}(b_{2, s})
\quad(r<s ) \;; \\
& & \; a_{2, s}^{b_{r}}={}^{[b_{2,r},
 \sigma_1^{-2}]}(a_{2, s})
\quad(r<s ) \;; \\
& & \; b_{2, s}^{b_{r}}={}^{[b_{2,r},
 \sigma_1^{-2}]}(b_{2, s})
\quad(r<s ) \;; \\
& & \;{}^{a_r}a_{2, s}={}^{[\sigma_1^2, a_{2, r}^{-1}]}(a_{2,s})\quad(s<r) \;; \\  
& & \;{}^{b_r}a_{2, s}={}^{[\sigma_1^2, b_{2, r}^{-1}]}(a_{2,s})\quad(s<r) \;; \\  
& & \;{}^{a_r}b_{2, s}={}^{[\sigma_1^2, a_{2, r}^{-1}]}(b_{2,s})\quad(s<r) \;; \\  
& & \;{}^{b_r}b_{2, s}={}^{[\sigma_1^2, b_{2, r}^{-1}]}(b_{2,s})\quad(s<r) \;; \\  
& &  b_{2, r}^{a_{r}}= (a_{2,r}\sigma_1^{-2}a_{2,r}^{-1}) b_{2, r}[\sigma_1^{-2}, a_{2,r}]
 \quad(1 \le r \le g ) \, ; \\
& & {}^{a_r}b_{2, r}=\sigma_1^2 b_{2, r} [a_{2, r}^{-1},\sigma_1^2]
 \quad(1 \le r \le g ) \,  ;\\
& & {}^{b_r}a_{2, r}=  a_{2,r} \sigma_1^{-2} 
 \quad(1 \le r \le g  ) \, ; \\
& &  a_{2, r}^{b_r}=  a_{2,r} b_{2, r} \sigma_1^2  b_{2, r}^{-1}
 \quad(1 \le r \le g  ) \, .
\end{eqnarray*}
A consequence of these identities and relation $(R1)$ is  that
 $\omega_i^{a_r}, {}^{a_r}\omega_i, \omega_i^{b_r}, {}^{b_r}\omega_i
\in \langle \mathcal{G} \rangle$ ($i, r=1, \dots,g$). 
\edim

\begin{lem} \label{sect2:bord}
Set also $\{\omega_1, \dots, \omega_{2g}, \tau_1, \dots \tau_{n-1}\}$ in
 $B^0(n, F)$ for
 $\{\phi_n(\omega_1), \dots,$
$ \phi_n(\omega_{2g}),$
$ \phi_n(\tau_1),
\dots \phi_n(\tau_{n-1})\}$. Then   $Ker(\theta)$ is  freely generated by
$\{\omega_1, \dots, \omega_{2g},$
$ \tau_1, \dots, $ $ \tau_{n-1}\}$.
\end{lem}
\pre
The diagram
\begin{diagram}
   P(n,F) &  \rTo^{\theta} & P(n-1,F) \\
   \dInto &            & \dInto \\
    B^0(n, F) &    \rTo^{\theta} &  B(n-1, F)
   \end{diagram} 
is commutative and the kernels of horizontal maps are the same.
As stated in section \ref{section},  $Ker(\theta)=
\pi_1(F \setminus \{ P_1, \dots, P_{n-1}\}, P_n)$.
If the fundamental domain is changed as in  Figure 6 and 
 $\omega_j, \tau_i$ are considered as loops of the fundamental group
 of $F \setminus \{P_1, \dots, P_{n-1} \}$ based on $P_n$, it is clear that  
$\pi_1(F \setminus \{P_1, \dots, P_{n-1}\}, P_n)=\langle
\omega_1, \dots, \omega_{2g}, \tau_1, \dots \tau_{n-1}| \; \emptyset
\; \rangle  $.
\edim

\begin{figure}[ht]
 \hspace{20pt}\psfig{figure=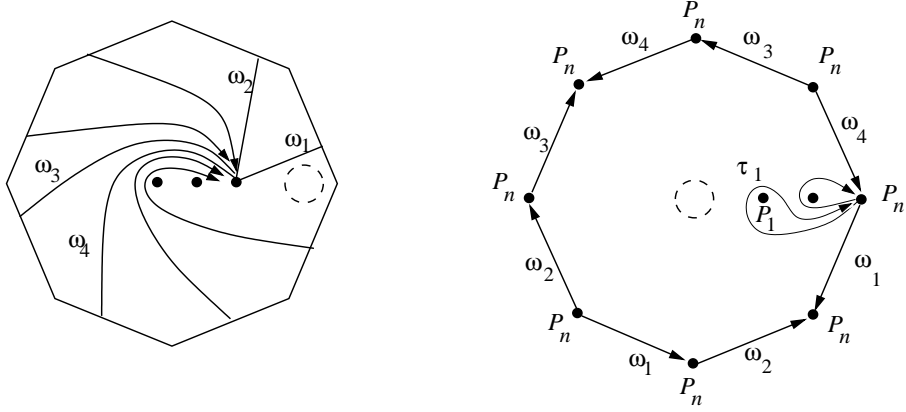,width=12cm} 
 \caption{Interpretation of $\omega_j, \tau_i$  as loops of the fundamental group.}
 \end{figure}

\begin{lem} \label{sect3:bord}
$\phi_{n|\widetilde{B}^0(n, F)}$ is an isomorphism.
\end{lem}
\pre
From the previous Lemmas  it follows that the map from
$Ker(\tilde{\theta})$
to $Ker(\theta)$ is an isomorphism. The Five Lemma and the inductive
assumption conclude the proof.
\edim


\subsection{End of the proof}

In order to show that $\phi_n$ is an isomorphism, let us remark
first that it is onto. In fact, from the previous Lemma the image of $\widetilde{B}(n,F)$ 
contains $P_n$ and on the other hand   $\widetilde{B}(n,F)$ 
surjects on $\Sigma_n$.
Since
the index of $B^0(n, F)$ in $B(n, F)$ is $n$, it is
sufficient to show that $[ \widetilde{B}(n,F):\widetilde{B}^0(n,F)]=n$.
Consider the elements $\rho_j=\sg_j \cdots \sg_{n-1}$  (we set $\rho_n=1)$
in  $\widetilde{B}(n,F)$. We claim that 
$\bigcup_i \rho_i \widetilde{B}^0(n,F)= \widetilde{B}(n,F)$. We  only have to show  that
for any (positive or negative) generator $g$ of  $\widetilde{B}(n,F)$ and $i=1, \dots, n$ there
exists $j=1, \dots, n$ and  $x \in \widetilde{B}^0(n,F)$ such that
$$
g \rho_i = \rho_j x \; .
$$
If  $g$ is a classical braid, this result is well-known
(\cite{cho}). Other cases 
come almost directly  from  the definition of
$\omega_j$.
Thus  every element of $\widetilde{B}(n,F)$ can be written in the form
$\rho_i \widetilde{B}^0(n,F)$. Since $\rho_i^{-1} \rho_j \notin
\widetilde{B}^0(n,F)$
for $i\not= j$ we are  done.
\edim

The  previous proof holds also for $p>1$.
This time 
$\widetilde{B}^0(n, F)$ is the subgroup of $\widetilde{B}(n, F)$
generated  
by $a_1, \dots, a_{g}, b_1, \dots, b_{g},\sg_1, \dots, \sg_{n-2},$ $ \tau_1, \dots ,
\tau_{n-1},$
$
\omega_1, \dots , \omega_{2g},$ $ \zeta_1, \dots ,  \zeta_{p-1}$
where
$\tau_j, \,\omega_r$ are defined as above and $\zeta_j = $
$\sg_{n-1}^{-1} \cdots \sg_1^{-1}$ $  z_j
\sg_1, \cdots, \sg_{n-1}$.
\edim

\section{Proof of Theorem \ref{men}}



\subsection{About the  section}
The steps of  the proof are the same.
We set again $B^0(n, F)=\pi^{-1}(\Sigma_{n-1})$. This time 
$\widetilde{B}^0(n, F)$ is the subgroup of $\widetilde{B}(n, F)$
generated  
by $a_1, \dots, a_{g},  b_1, \dots, b_{g},$ $ \sg_1, \dots, \sg_{n-2},$ $ \tau_1, \dots ,
\tau_{n-1},$
$
\omega_1, \dots , \omega_{2g},$
where
$\tau_j, \,\omega_r$ are defined as above.
Remark  that $\tau_1 \in \langle \mathcal{G} \rangle$ since from (TR) relation,
the following relation
$$
\tau_1=[\omega_1,\omega_2^{-1}] \cdots [\omega_{2g-1},\omega_{2g}^{-1}] \tau_{n-1}^{-1} \cdots \tau_{2}^{-1} \, ,
$$
holds in $\widetilde{B}^0(n, F)$.
When $F$ is a closed surface the corresponding $\tilde{\theta}$ has no section
(see section \ref{section}). Nevertheless, we are able to prove the analogous of Lemma
\ref{sect1:bord} (see section \ref{main:lema}).

\begin{lem} \label{sect1:closed}
Let  $F$ be  a closed surface. Then $Ker(\tilde{\theta)}$ is  generated by
$\{\omega_1, \dots,$ $ \omega_{2g},$
$ \tau_2, \dots \tau_{n-1}\}$.
\end{lem} 

The following Lemma is analogous to Lemma \ref{sect2:bord}. 

\begin{lem} \label{sect2:closed}
Let  $F$ be a  closed surface and set also $\{\omega_1, \dots, \omega_{2g}, \tau_2, \dots \tau_{n-1}\}$ in
 $B^0(n, F)$ for 
 $\{\phi_n(\omega_1), \dots,$
$ \phi_n(\omega_{2g}), \phi_n(\tau_2),
\dots \phi_n(\tau_{n-1})\}$. 
$Ker(\theta)$ is  freely generated by
$\{\omega_1, \dots, \omega_{2g}, \tau_2, \dots \tau_{n-1}\}$.
\end{lem}

Let $\rho_j=\sg_j \cdots \sg_{n-1}$  (where $\rho_n=1)$.
We may conclude  by checking that  for any 
generator $g$ of  $\widetilde{B}(n,F)$ (or its inverse) and $i=1, \dots, n$ there
exists $j=1, \dots, n$ and  $x \in \widetilde{B}^0(n,F)$ such that
$$
g \rho_i = \rho_j x \; ,
$$
which is a sub-case of previous situation.
\edim
 

\subsection{Proof of Lemma \ref{sect1:closed}} \label{main:lema}
To conclude the proof of Theorem \ref{men}, we give the demonstration
of Lemma \ref{sect1:closed}. Let us begin with the following Lemma.

\begin{lem} \label{main:lemma}
Let  $F$ be a closed surface and $\mathcal{G}=\{ \tau_2, \dots , \tau_{n-1},
\omega_1, \dots , \omega_{2g} \}$.
The subgroup $\langle \mathcal{G} \rangle$ is normal in $\widetilde{B}^0(n, F)$
\end{lem}
\pre
It suffices to consider relations in Lemma \ref{sect1:bord}.
 Remark that from relations shown in Lemma \ref{sect1:bord}, it follows also that  
 the set 
$$
\{ \gamma\tau_j \gamma^{-1} |j=1, \dots n-1, \, \gamma \; \mbox{word over}  \; \{\omega_1^{\pm 1}, \dots,
\omega_{2g}^{\pm 1}\}\} \, ,
$$ 
 is a system of generators for $\langle \langle \tau_1, \dots, \tau_{n-1} \rangle
\rangle \equiv \langle \langle\tau_{n-1} \rangle
\rangle $.
\edim

  In order to prove Lemma \ref{sect1:closed}, let us consider  the following diagram
\begin{diagram}
Ker \tilde{\theta}  & \rInto_i    & \widetilde{B}^0 (n, F)               & \rTo^{\tilde{\theta}}     & \widetilde{B}(n-1, F) \\
 \dTo^{t_n}         &             & \dTo^{q_n}                           & \ruTo^{\tilde{\theta}'}   &                       \\
Ker \tilde{\theta}' & \rInto_{i'} & \widetilde{B}^0(n, F)/\langle \langle \tau_{n-1} \rangle \rangle &       &
   \end{diagram} 

 In this diagram  $q_n$ is the natural projection, $\tilde{\theta}'$ is
defined by 
$\tilde{\theta}' \circ  q_n=\tilde{\theta}$ and $t_n$ is defined by
$i' \circ t_n=q_n \circ i$. Since $t_n$ is well defined and onto we deduce
that $Ker(t_n)= \langle \langle \tau_{n-1} \rangle
\rangle$.
 Now, $\tilde{\theta}'$ does have a natural section $s:\widetilde{B}(n-1, F) \to
 \widetilde{B}^0(n, F)/\langle \langle \tau_{n-1} \rangle \rangle$
defined as $s(a_i)=[a_j]$, $s(b_i)=[b_j]$  and $s(\sg_j)=[\sg_j]$, where $[x]$ is a representative
 of $x \in \widetilde{B}^0(n, F)$ in $\widetilde{B}^0(n, F)/\langle \langle \tau_{n-1} \rangle \rangle$.
Thus, using the same argument as in Lemma \ref{sect1:bord}, we derive that  
$Ker(\tilde{\theta}')=\langle \langle \mathcal{K} \rangle \rangle$,
  where $\mathcal{K}=\{ [\omega_1], \dots,
 [\omega_{2g}],[\tau_2],\dots[\tau_{n-1}] \}$.  From Lemma \ref{main:lemma}
it follows that  $\langle \mathcal{K} \rangle=\langle \langle \mathcal{K} \rangle \rangle$.
Moreover, since $\tau_i \in \langle \langle
 \tau_{n-1} \rangle \rangle$ for $i=1, \dots, n-2$, $Ker(\tilde{\theta}')=\langle [\omega_1], \dots,
 [\omega_{2g}] \rangle$.

From the exact sequence
$$
1 \to \langle \langle \tau_{n-1} \rangle
\rangle \to   Ker(\tilde{\theta})   \to Ker(\tilde{\theta'}) \to 1
$$ 
it follows that  $\{ \omega_1, \dots, \omega_{2g} \}$ and
a system of generators for $\langle \langle \tau_{n-1} \rangle
\rangle$   form a system of
generators
for  $Ker(\tilde{\theta})$. From the remark in Lemma \ref{main:lemma}
it follows that  $Ker(\tilde{\theta})=\langle  \tau_2, \dots , \tau_{n-1},
\omega_1, \dots , \omega_{2g} \rangle$.
\edim

\section{Other presentations and remarks} \label{otherpres}
\subsection{Braids on $p$-punctured spheres}
We recall that 
the exact sequence   
$$
1 \rTo \pi_1(F \setminus \{P_1, \dots,
  P_{n-1}\}, P_n) \rTo P(n,F) \rTo^{\theta} P(n-1,F) \to 1
$$
holds also when $F=S^2$ (\cite{fa}).
Thus, previous arguments may  be repeated  in the case of the sphere, to obtain a new proof for the well-known 
presentation of braid groups on the sphere as quotients of classical braid groups.
On the other hand, when $F$ is $p$-punctured sphere we have the following result.

\begin{thm}  \label{sp}
Let   $F$  be a $p$-punctured sphere. The group  $B(n, F)$ admits
the following presentation:  

\begin{itemize}
\item Generators:
$\sg_1,\dots,\sg_{n-1}, z_1, \dots, z_{p-1} \,.$
\item Relations:
\begin{itemize}
\item  Braid relations, i.e.
\begin{eqnarray*}
 \sip \siip \sip &=&  \siip \sip \siip \, ;\\
\sip \sjp &=& \sjp \sip \quad   \mbox{for} \;| i-j | \ge 2  \, .
\end{eqnarray*}
\item Mixed relations:
\end{itemize}
\end{itemize}
\begin{eqnarray*}
&(R1)&\; z_j \sg_i= \sg_i z_j        \quad  (i \not=1, j=1, \dots, p-1 )\, ;\\ 
&(R2)&\;\sg_{1}^{-1} z_j \sg_{1} z_l=z_l\sg_{1}^{-1} z_j \sg_{1} \quad (j=1, \dots, p-1 , \; j<l )\, ;\\
&(R3)&\;\sg_{1}^{-1} z_j \sg_{1}^{-1} z_j=z_j\sg_{1}^{-1} z_j \sg_{1}^{-1}\quad (j=1, \dots, p-1) \, ;\\ 
\end{eqnarray*}
\end{thm}
We remark that this presentation coincides with Lambropoulou's presentation \cite{lam}.
\subsection{Braids on non-orientable surfaces}
On the other hand, previous arguments hold to show the following 
Theorems for non-orientable surfaces.

\begin{thm} \label{nonorop}
Let   $F$  be a 
non-orientable $p$-punctured surface of genus $g\ge 1$. The group  $B(n, F)$ admits
the following presentation:  

\begin{itemize}
\item Generators:
$\sg_1,\dots,\sg_{n-1}, a_1, \dots,a_{g}, z_1,\dots, z_{p-1} \,.$
\item Relations:
\begin{itemize}
\item  Braid relations, i.e.
\begin{eqnarray*}
 \sip \siip \sip &=&  \siip \sip \siip \, ;\\
\sip \sjp &=& \sjp \sip \quad   \mbox{for} \;| i-j | \ge 2  \, .
\end{eqnarray*}
\item Mixed relations:
\end{itemize}
\end{itemize}
\begin{eqnarray*}
&(R1) &\;a_r \sg_i=\sg_i  a_r      \quad  (1 \le r \le g;\; i\not= 1 )\, ; \\
&(R2) &\;\sigma_1^{-1} a_r \sigma_1^{-1} a_{r}= a_r \sigma_1^{-1} a_{r}\sigma_1 \quad  (1 \le r \le g) \, ; \\
&(R3)&\;  \sigma_1^{-1} a_{s} \sigma_1 a_r = a_r \sigma_1^{-1} a_{s} \sigma_1\quad (  s < r) \, ; \\
&(R4) &\; z_j\sg_i= \sg_i z_j        \quad  (i \not=n-1, j=1, \dots, p-1 )\, ;\\ 
&(R5)&\;  \sigma_1^{-1} z_i \sigma_1 a_r =a_r\sigma_1^{-1} z_i\sigma_1    \quad  (1 \le r \le g;\; i=1,
\dots, p-1; \; n>1 )\,; \\
&(R7)&\;\sg_{1}^{-1} z_j \sg_{1} z_l=z_l\sg_{1}^{-1} z_j \sg_{1} \quad (j=1, \dots, p-1 , \; j<l )\, ;\\
&(R8)&\;\sg_{1}^{-1} z_j \sg_{1}^{-1} z_j=z_j\sg_{1}^{-1}  z_j \sg_{1}^{-1} \quad (j=1, \dots, p-1)  \,    .
\end{eqnarray*}
\end{thm}

\begin{thm} \label{nonoro}
Let  $F$ be  a closed non-orientable surface of genus $g\ge 2$. The group $B(n,F)$ admits the following presentation: 
\begin{itemize}
\item Generators:
$\sg_1,\dots,\sg_{n-1}, a_1, \dots,a_{g} \,.$
\item Relations:
\begin{itemize}
\item Braid relations as in Theorem \ref{my:th}.

\item Mixed relations:
\end{itemize}
\end{itemize}
\begin{eqnarray*}
&(R1) &\;a_r \sg_i=\sg_i  a_r      \quad  (1 \le r \le g;\; i\not= 1 )\, ; \\
&(R2) &\;\sigma_1^{-1} a_r \sigma_1^{-1} a_{r}= a_r \sigma_1^{-1} a_{r}\sigma_1 \quad  (1 \le r \le g) \, ; \\
&(R3)&\;  \sigma_1^{-1} a_{s} \sigma_1 a_r = a_r \sigma_1^{-1} a_{s} \sigma_1\quad (  s < r) \, ; \\
&(R4) &\; z_j\sg_i= \sg_i z_j        \quad  (i \not=n-1, j=1, \dots, p-1 )\, ;\\ 
&(TR) &\; a_1\cdots a_{g}=\tre \, .
\end{eqnarray*}
\end{thm}

\vsp

We give only a geometric interpretation of generators (\cite{gm}).
To represent a braid in $F$ we consider  the surface as a polygon, 
this time of $2g$ sides as in Figure 9, and we make an additional cut: define the path
$e$ as in the left hand of the Figure  9 and cut the poligon along it. 
We get $F$ represented as in the right hand side of the same figure, where we can also see how we choose
the points $P_1, \dots, P_n$.
We show generators in Figure 10. Generators $\sg_j$ and $z_j$ are as above. For all $r=1, \dots, g$, the braid $a_r$ consists on the first string passing through the $r$-th wall, while the other strings are trivial paths.
Relations can be easily verified drawing corresponding braids. The  (TR) relation in Theorem \ref{nonoro}
is treated in \cite{gm}.
We remark that setting $g=1$ the previous Theorem provides a presentation for braid groups on the projective plane
(see also \cite{vb}).

\begin{figure}[ht]
\hspace{20pt}\psfig{figure=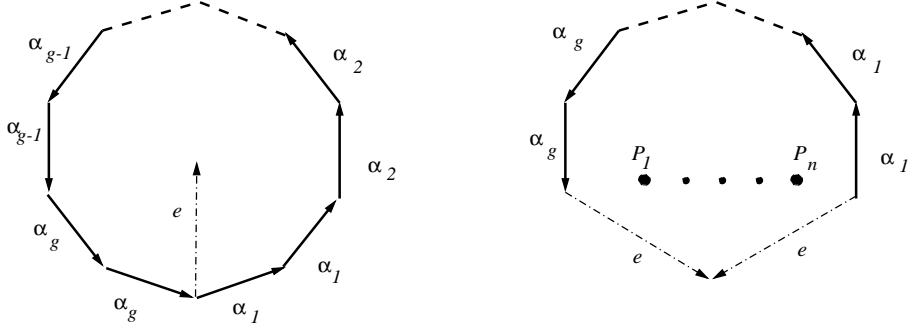,width=12cm} 
 \caption{Representation of a non-orientable surface $F$.}
 \end{figure}

\begin{figure}[ht]
\hspace{5pt}\psfig{figure=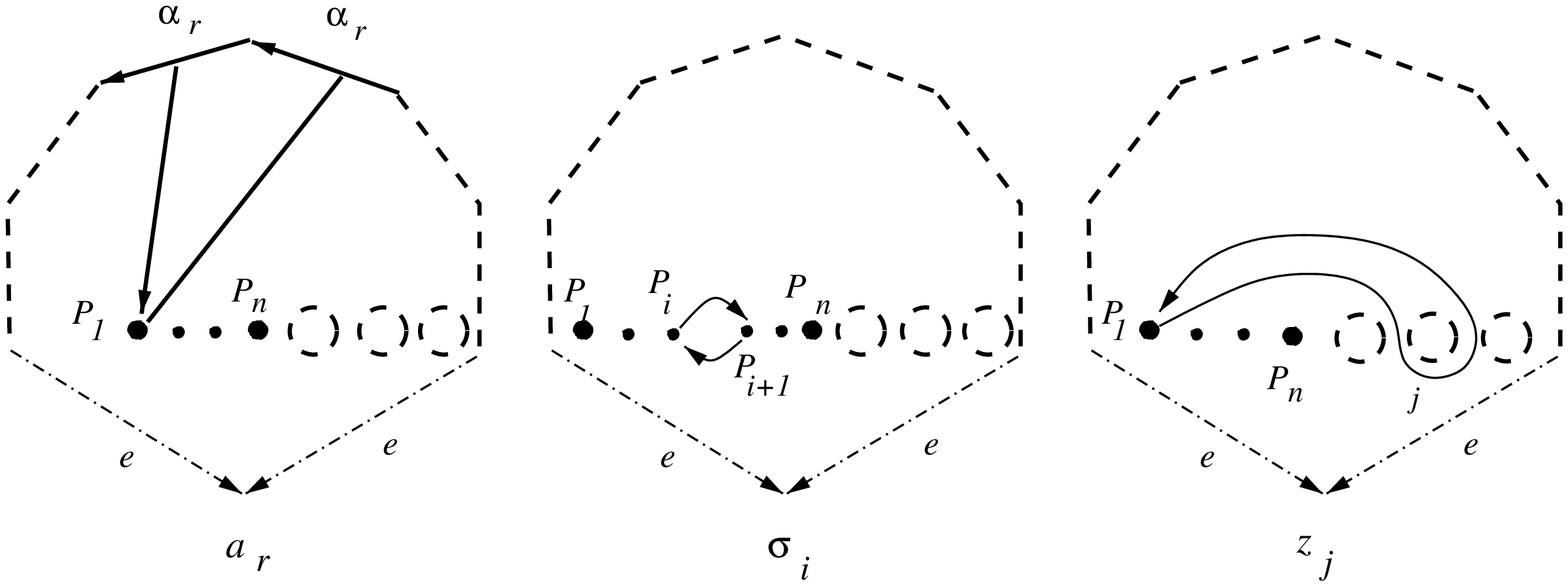,width=13cm} 
 \caption{Generators  as braids (for $F$ a non-orientable surface).}
 \end{figure}


\subsection{Gonz\'alez-Meneses' presentations} \label{sec}
Let  $F$ be  a closed orientable surface of genus $g\ge1$.
Using the same arguments outlined in previous sections we may provide an other presentation for $B(n,F)$.

\begin{thm}\label{meno}
Let  $F$ be  a closed orientable surface of genus $g\ge1$. The group $B(n,F)$ admits the following presentation: 
\begin{itemize}
\item Generators:
$\sg_1,\dots,\sg_{n-1}, b_1, \dots,b_{2g} \, .$
\item Relations:
\begin{itemize}
\item Braid relations as in Theorem \ref{my:th}.

\item Mixed relations:
\end{itemize}
\end{itemize}
\begin{eqnarray*}
&(R1) &\; b_r\sg_i= \sg_i b_r \quad (1 \le r \le 2g; i\not= 1 )\, ;\\
& (R2)&\;b_s \sg_1^{-1} b_{r}\sg_1^{-1} = \sg_1 b_{r} \sg_1^{-1} b_s \quad (1 \le s<r \le
2g)\, ; \\
&(R3) &\; b_{r} \sg_1^{-1} b_r \sg_1^{-1}=  \sg_1^{-1} b_r \sg_1^{-1} b_{r} \quad (1 \le r \le
2g)\, ;\\
&(TR) &\; b_1 b_{2}^{-1}\dots b_{2g-1} b_{2g}^{-1} b_1^{-1}  b_{2} \dots b_{2g-1}^{-1} b_{2g}=\tre \, .
\end{eqnarray*}
\end{thm}
A  closed orientable surface $F$ of genus $g\ge 1$
is represented as a polygon $L$ of $4g$ sides, where  
opposite edges are identified. Figure 1.10 gives a geometric interpretation
of generators. Relations can be  easily verified on corresponding braids.

\begin{figure}[ht]
\hspace{40pt}\psfig{figure=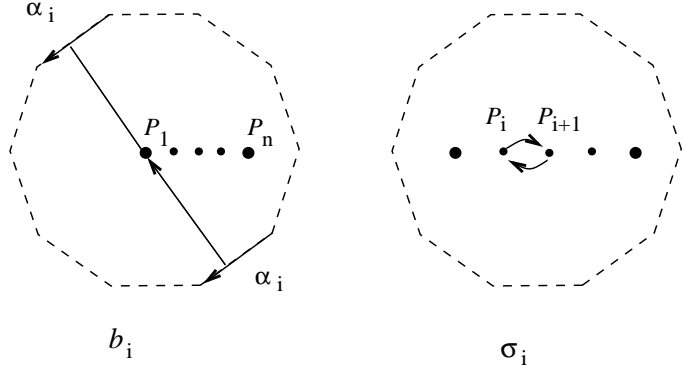,width=9cm} 
 \caption{Generators  as braids (for $F$ an orientable closed surface).}
 \end{figure}

The presentation  in Theorem \ref{meno} is close to 
Gonz\'alez-Meneses' presentation. 
\begin{thm}{\bf (\cite{gm})} \label{meneses}
Let  $F$ be  a closed orientable surface of genus $g\ge1$. The group $B(n,F)$ admits the following presentation: 
\begin{itemize}
\item Generators: $\sg_1,\dots,\sg_{n-1}, a_1, \dots,a_{2g} \, .$
\item Relations:
\end{itemize}
\begin{eqnarray*}
&(1)& \sip \siip \sip =  \siip \sip \siip \,, \\
&(2)& \sip \sjp = \sjp \sip \quad   \mbox{for} \;| i-j | \ge 2  \, , \\
&(3)& [ a_r, A_{2, s} ] = 1  \quad (1 \le r,s\le 2g; \; r \not= s)\, , \\
&(4)& [ a_r, \sg_i ] = 1 \quad (1 \le r \le 2g; i\not= 1 )\, ,\\
&(5)& [ a_1 \dots a_r, A_{2, r} ] = \sg_1^2 \quad (1 \le r \le 2g)\, ,\\
&(6)& a_1 \dots a_{2g} a_1^{-1}\dots a_{2g}^{-1}    = \tre \, ,
\end{eqnarray*}
where $A_{2, r}=\sg_1^{-1}(a_1\dots a_{r-1} a_{r+1}^{-1}\dots a_{2g}^{-1})\sg_1^{-1}\, .$
\end{thm}

Remark that the geometric interpretation of
$b_j$ corresponds to the braid generator $a_j$ when $j$ is odd and 
respectively to $a_j^{-1}$, when $j$ is even. 
Tedious 
computations show
that relations in Theorem \ref{meno} (after replacing generators 
$b_j$'s with $a_j$'s) imply  relations in Theorem \ref{meneses}.
In the same way, Theorem \ref{nonoro} can be also verified directly, checking
that  the relations in Theorem \ref{nonoro} imply all  relations  of the Gonz\'alez-Meneses' presentation 
for braid groups on non orientable closed surfaces in \cite{gm}.
However, we remark that the  presentation in Theorem \ref{nonoro} is simpler  and with less relations than  Gonz\'alez-Meneses' one.

On the other hand, it seems difficult to give an  algebraic proof of  the equivalence between presentation  in Theorem \ref{men} and
presentation in Theorem \ref{meneses}. 

\subsection{Applications}
We conclude this section with some remarks.
We recall that a subsurface $E$ of a surface $F$ is the closure of an open set of $F$. In order to avoid pathology, we  assume that
$E$ is connected and  that every boundary component of  $E$  either is a boundary component  of $F$ or lies in the interior of $F$.
We suppose also that $E$ contains $\mathcal{P}$. 
It is known \cite{rolpar} that  the natural map $\psi_n: B(n,E) \to B(n,F)$
induced by the inclusion
$E \subseteq F$ is injective if and only if $\overline{F\setminus E}$ does not contain a disk $D^2$. We may provide an analogous
characterisation about surjection. 

\begin{prop}
 Let $F$ be  a surface of genus $g\ge1$  with $p$ boundary components, and
 let $E$ be a subsurface of $F$. The natural map $\psi_n: B(n,E) \to B(n,F)$ induced by the inclusion
$E \subseteq F$ is surjective if and only if $\overline{F\setminus E} = \amalg D^2$.
 \end{prop}
\pre
Remark that 
the natural morphism 
$$
\psi_1: \pi_1(E,P_1) \rightarrow \pi_1(F,P_1)
$$
is a surjection if and only if  $\overline{F\setminus E} = \amalg D^2$.
Now consider a pure braid $p \in P(n,F)$ as a $n$-tuple of paths $(p_1, \dots, p_n)$ and
let $\chi: P(n,F) \to \pi_1(F)^n$  be the map defined by $\chi(p)=(p_1, \dots, p_n)$.
We   have the following commutative diagram 
\begin{diagram}
   P(n, E) & \rTo^{\chi} & \pi_1(E)^n\\
   \dTo^{\psi_n} &            & \dTo_{\psi_1 \times \cdots \times \psi_1} \\
  P(n, F) & \rTo^{\chi} & \pi_1(F)^n\\
\end{diagram} 
Since $\chi$ is surjective (\cite{bir2}) we deduce that $\psi_n$ is not surjective on  $ P(n, F)$  and thus 
on $B(n, F)$.
When $E$ is obtained from $F$ removing $k$ disks,  Theorems \ref{my:th}, \ref{men}, \ref{sp}
\ref{nonorop} and \ref{nonoro} give a description of 
$Ker(\psi_n)$. This result  can also be easily obtained 
from the remark that $B(n,E)$ is a  subgroup of $B(n+k,F)$ and that the map 
$\psi_n$ corresponds to the usual projection $B(n+k,F) \to B(n,F)$.
The existence of a braid combing in $B(n+k,F)$, analogous to that of the Artin braid group $B_n$,
implies the claim. 
\edim

\begin{prop} \label{hom}
Let $F$ be a orientable surface of genus $g \ge 1$, possibly with boundary.
Let $N_n(F)$ be the normal closure of   $B_n$ in  $B(n,F)$.
The quotient $B(n,F)/N_n(F)$ is isomorphic to $H_1(F)$, the first homology group of 
the surface $F$.
\end{prop}
\pre
Setting $\sg_j=1$ for $j=1, \dots, n-1$ in Theorems \ref{my:th} and \ref{men} we obtain
a presentation for $H_1(F)$.
\edim

\section{Surface pure braid groups }
Several presentations  for surface braid groups are known,
when $F$ is a closed surface or a holed disk (\cite{gm}, \cite{gua}, \cite{lam}, \cite{sco}).
In Theorem \ref{pure} we provide a  presentation for pure braid groups on orientable
surfaces with boundary. This presentation is   close to the 
standard presentation of the pure braid group $P_n$ on the disk.
We provide also the  analogous presentation for  pure braid groups on orientable
closed surfaces. 

Pure braid groups $P_n$
are a main ingredient in the construction of an universal finite type invariant 
for links in $\rea^3$ (see \cite{pap}). 

Using an approach similar to the case of $P_n$, in \cite{gonpar} 
it is constructed an universal finite type invariant for braids on orientable closed surfaces. 
 This construction is based 
on the group $K_n(F)$, the  normal closure of classical pure braid group $P_n$ in $P(n,F)$.

Consider the sub-surface $E$ obtained 
removing the  handles of $F$.
Let $Y_n(F)$ be the normal closure of $P(n,E)$ in $P(n,F)$.
Using classical techniques and  our presentation for surface pure braid groups, we prove 
that  the group $Y_n(F)$,  which contains properly $K_n(F)$,  is residually nilpotent (Theorem \ref{rtf}). 
On the other hand the group $Y_n(F)$ is the ``biggest subgroup'' of $P(n,F)$ on which one can use 
classical techniques, and the question whether $P(n,F)$ is residually nilpotent,
when  $F$ is a surface with genus, 
remains open.

\subsection{Presentations for surface pure braid groups }

\begin{thm} \label{pure}
Let  $F$ be an orientable surface
of genus $g\ge 1$ with $p>0$ boundary components. $P(n, F)$ admits the following presentation:
\begin{itemize}
\item Generators: 
$$\{A_{i,j}\; | \;1 \le i \le 2g + p + n -2,  2g +p\le j \le 2g + p + n -1, i<j \}.$$
\item Relations:
\end{itemize}
\begin{eqnarray*}
 &(PR1)&  A_{i,j}^{-1}  A_{r,s} A_{i,j} = A_{r,s} \;  \; \mbox{if} \, \,(i<j<r<s)  \;   \mbox{or} \, (r+1<i<j<s),\\ 
  &     &  \mbox{or} \, (i=r+1<j<s \, \, \mbox{for even} \, \, r<2g  \, \, \mbox{and}  \,
 \, r\ge 2g ) \,; \\ 
 &(PR2)&  A_{i,j}^{-1}  A_{j,s} A_{i,j} = A_{i,s}  A_{j,s} A_{i,s}^{-1} \;  \; \mbox{if} \, \, (i<j<s)\,;\\
 &(PR3)&  A_{i,j}^{-1}  A_{i,s} A_{i,j} = A_{i,s} A_{j,s} A_{i,s}  A_{j,s}^{-1} A_{i,s}^{-1} \; \;   \mbox{if} \, \, (i<j<s)\,; \\
 &(PR4)& A_{i,j}^{-1}A_{r,s}A_{i,j}=A_{i,s}A_{j,s}A_{i,s}^{-1}A_{j,s}^{-1}A_{r,s}A_{j,s}A_{i,s}A_{j,s}^{-1}A_{i,s}^{-1} \\
 &     &   \mbox{if} \, \,(i+1<r<j<s)  \;   \mbox{or} \\
 &     &   \, \, (i+1=r<j<s  \; \mbox{for odd} \, \, r<2g  \, \, \mbox{and}  \, \, r>2g )\,; \\
 &(ER1)&  A_{r+1,j}^{-1} A_{r,s} A_{r+1,j}=A_{r,s} A_{r+1,s} A_{j,s} A_{r+1,s}^{-1}  \\ 
 &     &     \mbox{if} \,  \,r \, \mbox{even and}\, \,r<2g                 \,               ; \\
 &(ER2)&   A_{r-1,j}^{-1}  A_{r,s} A_{r-1,j}=  A_{r-1,s} A_{j,s} A_{r-1,s}^{-1}  A_{r,s} A_{j,s} A_{r-1,s} A_{j,s}^{-1}A_{r-1,s}^{-1} \\ 
 &    &   \mbox{if} \,  \,r \, \mbox{odd and}\, \,r<2g \,.
\end{eqnarray*}
\end{thm}
\pre
The choice of the notation is motivated by the notation for standard generators of $P_n$ from \cite{bir}. 
Let $\widetilde{P}(n-1,F)$ be the group defined by above  presentation.
We give in Figure 1.11 a picture of  corresponding braids on the surface.
Let $h= 2g +p -1$. In respect of the presentation for  $B(n,F)$ given
in Theorem \ref{my:th}, the elements $A_{i,j}$ are the following braids:
\begin{itemize}
\item $A_{i,j}= \sg_{j-h} \cdots  \sg_{i+1-h} \sg_{i-h}^2 \sg_{i+1-h}^{-1}
 \cdots \sg_{j-h}^{-1}$, for $i \ge 2g + p$ ;
\item $A_{i,j}= \sg_{j} \cdots  \sg_{1} z_{i-2g}^{-1} \sg_{1}^{-1} \cdots \sg_{j-h}^{-1}$, for $2g < i < 2g + p$ ;
\item $A_{2i,j}= \sg_{j} \cdots  \sg_{1}  a_{g-i+1}^{-1} \sg_{1}^{-1} \cdots \sg_{j-h}^{-1}$, for 
$1 \le i \le g$ ;
\item $A_{2i-1,j}= \sg_{j-h} \cdots  \sg_{1}  b_{g-i+1}^{-1} \sg_{1}^{-1} \cdots \sg_{j-h}^{-1}$, for 
$1 \le i \le g$ .
\end{itemize}
The relations (PR1), \dots, (PR4) correspond to the
classical relations for $P_n$. The new relations arise when we consider two generators
$A_{2i,j}$, $A_{2i-1,k}$, for $1 \le i \le g$ and $j \not= k$. They correspond to two 
loops based at two different points which go around the same  
handle. Relations  (ER1)   and (ER2) can be  verified by explicit pictures or using relations in Theorem \ref{my:th}.
\begin{figure}[ht] 
\hspace{60pt}\psfig{figure=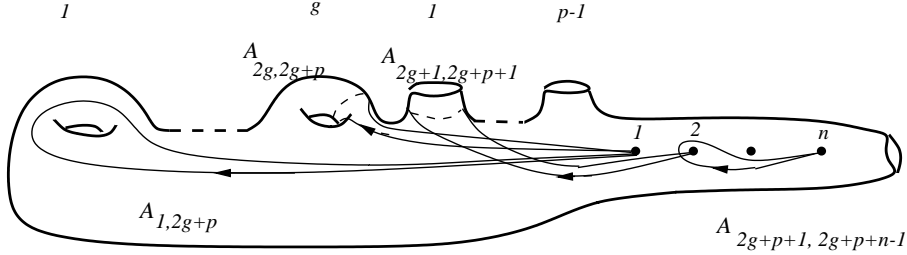,width=12cm} 
\caption{Geometric interpretation of $A_{i,j}$. We mark again with $A_{i,j}$ the only non trivial string of the braid
$A_{i,j}$}
\end{figure}
The technique to prove that $(PR1), \dots, (ER2)$ is a complete system of relations for $P(n,F)$
 is well known (\cite{gm}, \cite{gua}, \cite{lam}, \cite{sco}).
As shown in \cite{joh}, given an exact sequence
$$
1 \to A \to B \to C \to 1 \, ,
$$
and  presentations $\langle G_A, R_A \rangle$ and
$\langle G_C, R_C \rangle$, we can derive a presentation $\langle G_B, R_B \rangle $
for $B$, where $G_B$ is the set  of generators $G_A$ and  coset representatives of $G_C$.
The  relations  $R_B$ are given by  the union of three sets. The first corresponds to relations
$R_A$, and the second one  to writing each relation in $C$ in terms of corresponding coset representatives
as an element of $A$. The last set  corresponds to the fact that
the action under conjugation of each coset representative of generators of $C$ (and their inverses) on each generator
of $A$ is an element of $A$.
We can apply this result on (PBS) sequence.
 The presentation is correct for
$n=1$. By induction, suppose that for $n-1$, $\widetilde{P}(n-1,F) \cong P(n-1,F)$.
The set of elements
$A_{i,2g+n+p-1}$ $(i=1, \dots, 2g+n+p-2)$ is a system of generators for
$\pi_1(F \setminus \{P_1, \dots,
  P_{n-1}\}, P_n)$. 
To show that  $(PR1), \dots, (ER2)$ is a complete system of relations for $P(n,F)$ it
 suffices to prove that relations $R_{P(n,F)}$  are a consequence  of  relations
$(PR1), \dots, (ER2)$.
Since $\pi_1(F \setminus \{P_1, \dots,
  P_{n-1}\}, P_n)$ is a free group on the given generators, we just have to check the second and the third set
of relations. Consider as  coset representative for the generator $A_{i,j}$ in $P(n-1,F)$ the generator
$A_{i,j}$ in $P(n,F)$. Relations lift directly  to relations in $P(n,F)$.
The action of $A_{i,j}^{-1}$ on $\pi_1(F \setminus \{P_1, \dots,
  P_{n-1}\}, P_n)$ may be deduced from that of $A_{i,j}$.
In fact, relations (PR2) and (PR3) imply that
$$
 A_{i,j} A_{i,2g+n+p-1}  A_{j,2g+n+p-1}=  A_{i,2g+n+p-1}  A_{j,2g+n+p-1} A_{i,j} \; ,
$$
for all $i<j<2g+n+p-1$,
and from this relation and relations (PR2) we deduce that 
$$
A_{i,j} A_{i,2g+n+p-1} A_{i,j}^{-1}=A_{j,2g+n+p-1}^{-1} A_{i,2g+n+p-1} A_{j,2g+n+p-1} \; ,
$$
for all $i<j<2g+n+p-1)$.
It follows that
$$
A_{s,j} A_{i,2g+n+p-1} A_{s,j}^{-1} \in \langle A_{1,2g+n+p-1}, \dots,   A_{2g+n+p-2,2g+n+p-1} \rangle \, , $$
for all $s<j<2g+n+p-1$.

Thus we have proved that $ \langle A_{1,2g+n+p-1}, \dots,   A_{2g+n+p-2,2g+n+p-1} \rangle$
is a normal subgroup and that  also the third set of relations of $R_{P(n,F)}$ is a consequence of $(PR1), \dots,$
$ (ER2)$. 
\edim

In the same way we can prove the following Theorem.

\begin{thm} \label{pureclosed}
Let  $F$ be an orientable closed surface
of genus $g\ge 1$. $P(n, F)$ admits the following presentation:
\begin{itemize}
\item Generators: $\{A_{i,j}\; | \;1 \le i \le 2g+ n -1,  2g +1\le j \le 2g + n, i<j \}.$
\item Relations:
\end{itemize}
\begin{eqnarray*}
 &(PR1)&  A_{i,j}^{-1}  A_{r,s} A_{i,j} = A_{r,s} \;  \; \mbox{if} \, \,(i<j<r<s)  \;   \mbox{or} \, (r+1<i<j<s),\\ 
  &     &  \mbox{or} \, (i=r+1<j<s \, \, \mbox{for even} \, \, r<2g  \, \, \mbox{and}   \, \, r>2g  ) \,; \\
 &(PR2)&  A_{i,j}^{-1}  A_{j,s} A_{i,j} = A_{i,s}  A_{j,s} A_{i,s}^{-1} \;  \; \mbox{if} \, \, (i<j<s)\,;\\
 &(PR3)&  A_{i,j}^{-1}  A_{i,s} A_{i,j} = A_{i,s} A_{j,s} A_{i,s}  A_{j,s}^{-1} A_{i,s}^{-1} \; \;   \mbox{if} \, \, (i<j<s)\,; \\
 &(PR4)& A_{i,j}^{-1}A_{r,s}A_{i,j}=A_{i,s}A_{j,s}A_{i,s}^{-1}A_{j,s}^{-1}A_{r,s}A_{j,s}A_{i,s}A_{j,s}^{-1}A_{i,s}^{-1} \\
 &     &   \mbox{if} \, \,(i+1<r<j<s)  \;   \mbox{or} \\
 &     &   \, \, (i+1=r<j<s  \; \mbox{for odd }  \, \, r<2g  \, \, \mbox{and}   \, \, r>2g ) \,; \\ 
 &(ER1)&  A_{r+1,j}^{-1} A_{r,s} A_{r+1,j}=A_{r,s} A_{r+1,s} A_{j,s} A_{r+1,s}^{-1}  \\ 
 &     &     \mbox{if} \,  \,r \, \mbox{even and}\, \,r<2g                 \,               ; \\
 &(ER2)&   A_{r-1,j}^{-1}  A_{r,s} A_{r-1,j}=  A_{r-1,s} A_{j,s} A_{r-1,s}^{-1}  A_{r,s} A_{j,s} A_{r-1,s} A_{j,s}^{-1}A_{r-1,s}^{-1} \\ 
 &    &   \mbox{if} \,  \,r \, \mbox{odd and}\, \,r<2g   \,               ; \\
 &(TR)&    [ A_{2g, 2g+k}^{-1}, A_{2g-1, 2g+k}] \cdots [ A_{2, 2g+k}^{-1}, A_{1, 2g+k}] =
   \prod_{l=2g+1}^{2g+k-1} A_{l,2g+k} \times \\
 & &  \times \prod_{j=2g+k+1}^{2g+n} A_{2g+k,j}   \;  \; k=1, \dots, n    \,.
\end{eqnarray*}
\end{thm}

\begin{rem}\label{emb}
Let $E$ be a holed disk.
 Theorem \ref{pure} provides a presentation 
for $P(n,E)$ (\cite{lam}).
Let us recall that  $P(n,E)$ is a (proper) subgroup
of $P_{n+k}$, where $k$ is the number of holes in $E$.
\end{rem}

\begin{rem} \label{emb2}
We recall that $P_n$ embeds in $P(n,F)$ (\cite{rolpar})
and thus  $P_n$ is isomorphic to the subgroup 
$$
 \langle A_{i,j} \, | \, 2g + 1 \le i < j \le 2g  + n  \rangle \, , 
$$
when $F$ is a closed surface and 
$P_n$ is isomorphic to  
$$
P_n= \langle A_{i,j} \, | \, 2g + p \le i < j \le 2g  + p + n -1  \rangle \, , 
$$
when $F$ is a surface with $p>0$ boundary components.
Consider the sub-surface $E$ obtained 
removing $g$ handles from $F$.
The group    $P(n,E)$ embeds in $P(n,F)$ (\cite{rolpar})
and it is isomorphic to the subgroup 
$$
 \langle  A_{i,j} \cup A_{2k-1,l} \cup  A_{2k,l}^{-1} A_{2k-1,l}^{-1}  A_{2k,l}\, | \, 2g + 1 \le i < j \le 2g + n \; ,
\; 1 \le k \le g  \rangle  \, ,  
$$
when $F$ is a closed surface and respectively to the subgroup 
$$
 \langle  A_{i,j} \cup A_{2k-1,l} \cup  A_{2k,l}^{-1} A_{2k-1,l}^{-1}  A_{2k,l}\, | \, 2g+ 1 \le i  < j \le 2g + p + n-1 \; ,
\;  1 \le k \le g \rangle  \, , 
$$
when $F$ is a surface with $p>0$ boundary components.
\end{rem}

\begin{rem} 
When $F$ is a surface with genus,
from relation (ER1) we deduce that 
generators $A_{i,j} $ for $2g+p \le i <j \le 2g+n+p-1$,
which  generate  a subgroup isomorphic to $P_n$,
are redundant. Then  Theorem \ref{pure}
provides a (homogeneous) presentation for 
$P(n,F)$ with $(2g+p-1)n$ generators.
\end{rem}

\subsection{Remarks on the normal closure of  $P_n$ in $P(n,F)$}

As corollary of previous presentations we  give an easy proof of a well-known fact on
$K_n(F)$
(\cite{gol}).

\begin{lem}  \label{lem:fond} Let $\chi: P(n,F) \to \pi_1(F)^n$ be the map defined by $\chi(p)=(p_1, \dots, p_n)$.
Let $F$ be a closed  orientable surface possibly  with boundary.
Let $K_n(F)$ be  the normal closure of $P_n$ in $P(n,F)$. Then
$$Ker(\chi)=K_n(F) \, .$$
 \end{lem}
   \pre 
The set $\{\chi(A_{i,j}) \, | \,1 \le i \le 2g + p-1$ and  $2g + p  \le j \le 2g + p + n -1\}$  
forms a complete set of generators for  $
        \pi_1(F)^n$. On the other hand, from relation (ER1) it follows that
   $\chi(A_{i,j})=1$ for $2g + p  \le i< j  \le 2g + p + n -1$ and thus 
$$
Ker(\chi)= \langle \langle A_{i, j} \, | \, 2g + p \le i < j\le 2g + p + n -1 \rangle \rangle \, .
$$
Remark \ref{emb2} concludes the proof.
\edim

\begin{prop} 
Let $F$ be an orientable surface possibly with boundary.
When $F$ is a torus 
$$[ P(n,F), P(n,F)] = K_n(F) \, .$$
Otherwise the strict inclusion holds:
$$[ P(n,F), P(n,F)] \supset K_n(F) \, .$$
\end{prop}
\pre
The inclusion $K_n(F) \subset  [ P(n,F), P(n,F)]$
is clear.

\noindent Suppose that $[ P(n,F), P(n,F)]=K_n(F)=Ker(\chi)$ for $g>1$.
It follows that  $\displaystyle \frac{P(n,F)}{Ker(\chi)}$ is abelian. This is false since $\pi_1(F)^n$ is not abelian
for $g>1$. 
Let $w \in [ P(n,F), P(n,F)]$. The sum of exponents $A_{i,j}$ in $w$ must be zero. 
The projection of $\chi(w)$ on any coordinate is the sub-word  of $w$ consisting of the generators associated to corresponding 
strand. Since the sum of exponents is zero, if $F$ is a torus this projection is trivial and the claim follows.
\edim

\noindent We recall that $Y_n(F)$ is the normal closure of $P(n,E)$ in $P(n,F)$, where $E$ is  the sub-surface
obtained 
removing the  handles of $F$.

\begin{prop} 
Let $F$ be an orientable surface with $p>0$ boundary components.
The following inclusions hold:
$$
K_n(F) \subset Y_n(F)\subset K_{n+2g+p-1}(F) \, .
$$
\end{prop}
\pre
 Remark \ref{emb2} and Lemma \ref{lem:fond} imply that the inclusion $K_n(F) \subset Y_n(F)$ is proper.
Since $P(n,E)$ is isomorphic to a subgroup of $P_{n+2g+p-1}$
and  $P(n,F)$ embeds in  $P(n+2g+p-1, F)$, it follows that $Y_n(F)$
is isomorphic to  a subgroup of $ K_{n+2g+p-1}(F)$.
\edim

\subsection{Almost-direct products} \label{residual}

It is known that 
\begin{itemize}
\item $\bigcap_{d=0}^{\infty}I(P_n)^d=\{0\}$; 
\item $I(P_n)^d/I(P_n)^{d+1}$ is a free $\z$-module for all $d\ge0$,
where  $I^k$
means the $k$-th power of the augmentation ideal of the group ring of
$P_n$.
\end{itemize}
This result follows from a more general statement 
on \emph{almost-direct products}
(see \cite{gonpar} or \cite{pap} ).

\begin{df}
Let $A, C$ be two groups. If   $C$ acts on $A$ and the induced action on the abelianization of $A$ is trivial,
we say that $A \rtimes C$ is an \emph{almost-direct product} of $A$ and $C$.
\end{df}

\begin{prop} \label{papadima}
Let $A, C$ be two groups. If   $C$ acts on $A$ and the induced action on the abelianization of $A$ is trivial,
$$
I(A \rtimes C)^m = \sum_{k=0}^{m} I(A)^k \otimes I(C)^{m-k} \quad \mbox{for all} \;  \;m \ge 0 \, .
$$ 
Let $B$ be  a finitely iterated almost-direct product of  free groups, then
\begin{itemize}
\item $\bigcap_{d=0}^{\infty}I(B)^d=\{0\}$; 
\item $I(B)^d/I(B)^{d+1}$ is a free $\z$-module for all $d\ge0$.
\end{itemize}
\end{prop}
As consequence of Theorem \ref{pure}
we  can prove that:  
\begin{prop} \label{rtf}
Let $F$  be an orientable surface  with boundary.
Then 
\begin{itemize}
\item $\bigcap_{d=0}^{\infty}I( Y_n(F))^d=\{0\}$; 
\item $I( Y_n(F))^d/I( Y_n(F))^{d+1}$ is a free $\z$-module for all $d\ge0$.
\end{itemize}
\end{prop}
\pre
We sketch a proof for $F$  orientable surface  with one boundary component.
Let $\pi_1(F)^n$ be provided  with presentation
\begin{eqnarray*}
\langle A_{j, 2g+k}\, \, j=1, \dots, 2g, \, k=1, \dots, n |  [ A_{j, 2g+k}, A_{l, 2g+q}] =1 \; \mbox{for all} \; 
j, l=1, \dots, 2g, \, 1 \le k \not=q \le n \, \rangle \, ,
  \end{eqnarray*}
where $A_{j, 2g+k}$  are the loops defined in Theorem \ref{pure}.
Let $F_{g,n}$ be the group with presentation 
\begin{eqnarray*}
\langle \widetilde{A_{j, 2g+k}}\, \,j=1, \dots, g, \, k=1, \dots, n | [  \widetilde{A_{j, 2g+k}},  \widetilde{A_{l, 2g+q}}] =1 \;
 \mbox{for all} \; j, l=1, \dots, g, \, 1 \le k \not=q \le n \, \rangle \, .
 \end{eqnarray*}
Let $\mu: \pi_1(F)^n  \to F_{g,n}$ be the map defined 
by $\mu(A_{2i-1, 2g+k})=\widetilde{A_{i, 2g+k}}$ and $\mu(A_{2i, 2g+k})=1$. 
One can proceed as in Lemma \ref{lem:fond} for showing that
$Ker(\mu \circ \chi)=Y_n(F)  \, $. Thus the following commutative diagram holds:

 \begin{diagram}
 &      &   1        &     &  1                 &              & 1                       &      &\\
 &      &   \uTo     &     &     \uTo           &              &  \uTo                   &      &       \\
1&\rTo  & F_{g,1}     &\rTo &  F_{g,n}  & \rTo         &  F_{g,n-1}    &\rTo  &1\\
 &      &   \uTo     &     &    \uTo ^{\mu \circ \chi}            &              &       \uTo^{\mu \circ \chi}               &      &   \\
1&\rTo  & Ker(\theta)&\rTo & P(n, F)            & \rTo^{\theta}& P(n-1, F)               & \rTo &1\\
 &      & \uTo       &     &     \uTo           &              &  \uTo                   &      &   \\
1&\rTo  &  G_n       &\rTo &  Y_n(F)            &\rTo^{\theta} &  Y_{n-1}(F)             &\rTo  &1 \\
 &      & \uTo       &     &      \uTo          &              &  \uTo                   &      &   \\
 &      &    1       &     &      1             &              &    1                    &      &     
    \end{diagram}
where $F_{g,1} $ is the free group on $g$ generators and
$G_n= Y_n(F) \cap Ker(\theta)$ is a free group.

\begin{lem} \label{gen}
  The following set is a  system of generators for $G_n$.
$$\{\gamma A_{j,2g+n} \gamma^{-1} | 2g < j <  2g+n\,
 \; \, \mbox{and} \; \,   1 \le j < 2g , \, j \; \, \mbox{even} \; \,   \}
$$
where $\gamma$ is a word on $\{A_{2k-1,2g+n}^{\pm{1}} |  1 \le k \le g   \}$.
\end{lem} 
\pre 
Consider the vertical sequence 
\begin{diagram}
1 &\rTo  &  G_n  & \rTo  &  Ker(\theta) &\rTo   & F_{g,1}     &\rTo &1 \, . 
\end{diagram}
Recall that $Ker(\theta) = \pi_1(F \setminus \{P_1, \dots, P_{n-1}\}, P_n)$.
A set of free generators for this group is given by
$\{A_{j,2g+n} |  1 \le j < 2g+n \}$.
The map  $\mu \circ \chi$ sends $A_{j,2g+n}$ in $1$ for $2g < j <  2g+n$
and  $1 \le j < 2g, \, j$ even. On the other hand,
$\mu \circ \chi(A_{2k-1,2g+n})=\widetilde{A_{k, 2g+1}}$ for $k=1, \dots, g$. 
\edim

Recall that the existence of a section for $\theta$ implies that  $Y_{n-1}(F)$ acts by conjugation on $G_n$ and thus on the 
abelianization $G_n/[G_n,G_n]$. The following Lemma and Proposition \ref{papadima} conclude the proof.
\edim

\begin{lem} \label{conj}
  The action of $Y_{n-1}(F)$ by conjugation on $G_n/[G_n,G_n]$ is trivial. 
\end{lem} 
\pre Let
$t \in \{A_{j,k} |  2g < k <  2g+n, 2g < j <  k \,
 \; \, \mbox{and} \; \,   1 \le j < 2g   \,, \;\mbox{j even}  \,   \}$
and $f \in  \{ A_{j,2g+n} |  2g < j <  2g+n\,
 \; \, \mbox{and} \; \,   1 \le j < 2g   \,, \;\mbox{j even}  \,   \}$. 
We need 
to verify that every $t$ acts trivially on $G_n / [G_n,G_n]$.
Presentation in Theorem \ref{pure} shows that
$$
(A) \quad t f t^{-1} \equiv f \quad (\mbox{mod} \, [G_n,G_n] ) \, ,$$
for every $t$ and $f$.
Now consider the action of  $t$ on $A_{2s,2g+n}$, for $s=1, \dots, g$.
We refer once again to  Theorem \ref{pure} 
for showing   that for every $t \in \{A_{j,k} |  2g < k <  2g+n, 2g < j <  k \,
 \; \, \mbox{and} \; \,   1 \le j < 2g , \, j \; \, \mbox{even}  \,   \}$,
 $$
 (B) \quad t  A_{2s-1,2g+n} t^{-1} =  h  A_{2s-1,2g+n} \quad ( 1 \le s\le g) \, , $$
 where $h \in G_n$.
Let $\gamma$ be a word on $\{A_{2k-1,2g+n}^{\pm{1}} |  1 \le k < g   \}$.
 From (A) and (B) it follows that, for every 
$t \in \{A_{j,k} |  2g < k <  2g+n, 2g < j <  k \,
 \; \, \mbox{and} \; \,   1 \le j < 2g , \, j  \; \, \mbox{even}  \,   \}$,
$ t \gamma f  \gamma^{-1} t^{-1}= t \gamma t^{-1} t f t^{-1} t \gamma^{-1} t^{-1}=
h \gamma  t f t^{-1} \gamma^{-1}  h^{-1}  \equiv 
h \gamma f  \gamma^{-1}  h^{-1}\equiv \gamma f  \gamma^{-1} \,$, where $h$ is an element of
$G_n$.
\edim

\begin{rem}
We notice that classical 
techniques do not apply to the whole group  $P(n,F)$.
The main obstruction is that, even when  the exact sequence (PBS)
splits, the action of  $P(n,F)$ on 
the abelianisation of $\pi_1(F \setminus \{x_1, \dots, x_{n-1}\})$ is not  trivial, because  of relations (ER1)
and (ER2). In particular,  when $F$ is a surface of genus $g \ge 1$,
it is presently unknown whether the graded group associated to the 
lower central series of $P(n,F)$ is torsion free. 
 \end{rem}



\newpage


\begin{thebibliography}{99} 
 
 \bibitem{bir}
J.~S.~Birman,
\newblock Braids, Links, and Mapping Class Groups,
\newblock Ann.  Math. Stud., Princeton Univ. Press, vol.82,  1974. 
 
 \bibitem{bir2}
J.~S.~Birman,
\newblock On Braids Groups,
\newblock {\em Comm. Pure and App. Math.}, 22 (1969), 213-238.

 \bibitem{cho}
W.~L.Chow,
\newblock On the algebraical braid group,
\newblock {\em Ann. Math.}, 49 (1948),  654-658.
 
 
\bibitem{fan}
E.Fadell and L.Neuwirth,
\newblock Configuration spaces,
\newblock {\em Math. Scand.}, 10 (1962), 111-118.

 


 \bibitem{fa}
E.Fadell and J.Van Buskirk,
\newblock The braid groups of $E^2$ and $S^2$,
\newblock {\em Duke Math. J.}, 29 (1962), 243-258.

\bibitem{falk}
M. Falk and R. Randell,
\newblock Pure braid groups and products of free groups,
 \newblock Contemporary Mathematics, 78 (1988), 217-228.

\bibitem{gol}
C.H. Goldberg,
\newblock An exact sequence of braid groups, 
\newblock  {\em Math. Scand.},  33 (1973), 69-82.

\bibitem{gm} 
J. Gonz\'alez-Meneses, 
 \newblock New presentations of surface braid groups,
\newblock  {\em J. Knot Theory Ramifications},  to appear, math.GT/9910020.  
  
\bibitem{gon} 
J. Gonz\'alez-Meneses,
\newblock Ordering pure braid groups on closed surfaces,
\newblock {\em preprint}, math.GT/0006155.

\bibitem{gonpar}
J.Gonz\'alez-Meneses and L.Paris,
\newblock Vassiliev Invariants for braids on surfaces,
\newblock {\em Trans. A.M.S.}, to appear, math.GT/0006014.

\bibitem{gua}
J. Guaschi and  P. Gon\c{c}alves,
\newblock On the structure of surface pure braid groups,
\newblock {\em Preprint}.

\bibitem{joh}
D.L. Johnson,
\newblock Presentation of groups,
\newblock LMS Lectures Notes 22 (1976), Cambridge University Press.

\bibitem{kul}
V.S. Kulikov and I. Shimada, 
\newblock On The Fundamental Groups of Complements to Dual Hypersurfaces of Projective Curves
\newblock  {\em  Preprint of the Max-Planck-Institut für Mathematik},  32 (1996).

\bibitem{lam} S.Lambropoulou, 
\newblock Braid structures in knot complements, handlebodies and 3-manifolds.
\newblock Knots in Hellas '98 (Delphi),  Ser.
Knots Everything, 24(2000),  274-289.

\bibitem{mor}
J.Morita,
\newblock A Combinatorial proof for Artin's Presentation of the
braid group and some cyclic analogue,
\newblock {\em Tsukuba J. Math.}, 16 (1992), no. 2, 439-442.

\bibitem{pap}
\c{S}. Papadima, 
\newblock The universal finite type invariant for braids, with integer coefficients,
\newblock {\em Topology and its Appl.}, \textbf{118} (2002), 169-185.

\bibitem{rolpar} 
L.Paris and D.~Rolfsen,
\newblock Geometric Subgroups of surface Braid Groups,
\newblock  {\em Ann.Inst. Fourier}, 42 (1998), 417-472.


\bibitem{sco}
G.P.Scott,
\newblock Braids groups and the group of homeomorphisms of
a surface,
\newblock {\em Proc.Camb. Phil. Soc.}, 68 (1970) 605-617.

 \bibitem{vlo}
V.~ Sergiescu,
 \newblock A direct approach to the planar graph presentations of the braid group, 
 \newblock {\em Singularities of Holomorphic Vector Fields and Related Topics}, Proceedings RIMS Kyoto,
878(1994), 103-107. 


 \bibitem{vb}
J.~ Van Buskirk,
 \newblock Braid groups of compact $2$-manifolds with elements of finite order, 
 \newblock {\em Trans. Amer. Math. Soc.}, 122 (1966), 81-97. 
 
 
\end{thebibliography}
\end{document}